# Continued fractions related to Narayana polynomials

Johann Cigler

**Abstract**

The generating functions of some sequences of Catalan numbers and Narayana polynomials have simple expansions as continued fractions of Jacobi type. We give an overview of these facts and prove analogous results for q-Narayana polynomials at q=-1.

## 1. Introduction

Let $(a_n)_{n\geq 0}$ be the sequence of Catalan numbers or the sequence of Narayana polynomials.

The generating functions of the sequences $(a_n)_{n\geq 0}$, $(a_{n+1})_{n\geq 0}$ and $\left(1,0,a_1,0,a_2,0,a_3,0,\cdots\right)$ have simple expansions as continued fractions of Jacobi type

(1)
$$\cfrac{1}{1-s_0z-\cfrac{t_0z^2}{1-s_1z-\cfrac{t_1z^2}{1-\ddots}}}.$$

We give an overview of these facts and prove analogous results for q-Narayana polynomials at q=-1.

## 2. Motivation

In order to motivate the consideration of these continued fractions let us sketch the context where they naturally appear.

Let $\left(a_n\right)_{n\geq 0}$ be a sequence of real numbers with $a_0=1$ and define a linear functional $L$ on $\mathbb{R}[x]$ by $L\left(x^n\right)=a_n$. Then the following result holds (cf. e.g.[3], [4],[5]):

If $\det\left(a_{i+j}\right)_{i,j=0}^{n-1}\neq 0$ for all $n\in\mathbb{N}$ there exists a sequence of monic polynomials $p_n(x)$ which are orthogonal with respect to $L.$ By a theorem of Favard there exist sequences $\left(s_n\right)$ and $\left(t_n\right)$ such that

(2)
$$p_n(x)=(x-s_{n-1})p_{n-1}(x)-t_{n-2}p_{n-2}(x).$$

If we define $a(n,k)$ by

(3)
$$a(n,k)=a(n-1,k-1)+s_ka(n-1,k)+t_ka(n-1,k+1)$$

with $a(n,k)=0$ for $k<0$ and $a(0,k)=\left[k=0\right]$ then $a(n,0)=a_n$ and the generating function $\sum_{n\geq 0}a_nz^n$ has a continued fraction expansion of the form (1).

The numbers $a(n,k)$ can be interpreted as weight of the set of Motzkin paths from $(0,0)$ to $(n,k).$

Recall that a Motzkin path is a lattice path in $\mathbb{Z}^2$ starting at $(0,0)$ consisting of up-steps $U=(1,1)$, horizontal steps $H=(1,0)$ and down-steps $D=(1,-1)$ which never goes below the $x-$axis. For given sequences $(s_k)$ and $(t_k)$ we define a weight $w$ on the Motzkin paths by $w(U)=1$ for up-steps $U$, $w(H)=s_k$ for horizontal steps $H$ on height $y=k$ and $w(D)=t_k$ for down-steps $D:(i-1,k+1)\to(i,k)$. The weight of a path is the product of the weights of its steps and the weight of a set of paths is the sum of the weight of its paths.. A Motzkin path from $(0,0)$ to $(2n,0)$ without horizontal steps is called a Dyck path of semi-length $n$.

## 3. Sequences of Catalan numbers

Let us start with the well-known situation for the Catalan numbers $C_n=\frac{1}{n+1}\binom{2n}{n}$.

They have the generating function

(4) $$C(z)=\sum_{n\geq 0}C_n z^n=\frac{1-\sqrt{1-4z}}{2z},$$

which satisfies

(5) $$C(z)=1+zC(z)^2$$

and implies

(6) $$C\left(z^2\right)=\frac{1}{1-z^2C\left(z^2\right)},$$

which gives

**Example 1**

*The generating function*

(7) $$C\left(z^2\right)=\sum_{n\geq 0}C_n z^{2n}$$

*has a continued fraction expansion with*

(8) $$s_n=0 \text{ and } t_n=1,$$

*This is equivalent with the well-known fact that* $C_n$ *can be interpreted as the number of Dyck paths of semi-length* $n$*.*

**Example 2**

*For*

(9) $$G(z)=\sum_{n\geq 0}C_{n+1}z^n=C(z)^2$$

*we get*

(10) $$s_n = 2, \quad t_n = 1.$$

*This also implies Touchard's identity*

(11) $$C_{n+1} = \sum_{k=0}^{\left\lfloor \frac{n}{2} \right\rfloor} \binom{n}{2k} 2^{n-2k} C_k.$$

**Proof**

From

$$G(z) = C(z)^2 = \frac{C(z)-1}{z} = \sum_{n\geq 0} C_{n+1} z^n \text{ we get } G(z) = C(z)^2 = \left(1 + zC(z)^2\right)^2 = (1 + zG(z))^2$$

and thus

(12) $$G(z)\left(1 - 2z - z^2 G(z)\right) = 1$$

which gives the continued fraction

(13) $$G(z) = C(z)^2 = \cfrac{1}{1 - 2z - \cfrac{z^2}{1 - 2z - \cfrac{z^2}{\ddots}}}.$$

(13) shows that $C_{n+1}$ can be interpreted as the weight of the Motzkin paths from $(0,0)$ to $(n,0)$ where the horizontal steps have weight $w(H) = 2$ and the down-steps have weight $W(D) = 1.$

To obtain the weight of the paths with $k$ up- and down-steps we can choose their positions in $\binom{n}{2k}$ ways and order them in $C_k$ ways. This gives Touchard's identity (11).

From $C(z) = 1 + zC(z)^2$ we get

(14) $$C(z) = \frac{1}{1 - zC(z)} = \frac{1}{1 - z\left(1 + zC(z)^2\right)} = \frac{1}{1 - z - z^2 G(z)}.$$

This implies

**Example 3**

*For* $C(z)$ *we get* $s_n = 2$ *for* $n > 0$, $s_0 = 1$ *and* $t_n = 1.$

## 4. Sequences of Narayana polynomials

The Narayana polynomials

(15) $$C_n(t)=\sum_k \binom{n}{k}\binom{n-1}{k}\frac{t^k}{k+1}$$

can be considered as refinements of the Catalan numbers $C_n(1)=C_n$.

Their generating function

(16) $$C(t,z)=\sum_{n\geq 0} C_n(t)z^n$$

satisfies (cf. [6])

(17) $$C(t,z)=\frac{1+(t-1)z-\sqrt{\left(1+(t-1)z\right)^2-4tz}}{2tz}=\frac{1}{1+(t-1)z}\frac{1-\sqrt{1-\dfrac{4tz}{\left(1+(t-1)z\right)^2}}}{\dfrac{2tz}{\left(1+(t-1)z\right)^2}},$$

i.e., in terms of $C(z)$

(18) $$C(t,z)=\frac{1}{1+(t-1)z}C\left(\frac{tz}{\left(1+(t-1)z\right)^2}\right).$$

By (5) this implies

$$\left(1+(t-1)z\right)C(t,z)=C\left(\frac{tz}{\left(1+(t-1)z\right)^2}\right)=1+\frac{tz}{\left(1+(t-1)z\right)^2}C\left(\frac{tz}{\left(1+(t-1)z\right)^2}\right)^2=1+tzC(t,z)^2$$

and

(19) $$C(t,z)=1+(1-t)zC(t,z)+tzC(t,z)^2.$$

From

(20) $$C\left(t^2,z^2\right)=1+(1-t^2)z^2C\left(t^2,z^2\right)+t^2z^2C\left(t^2,z^2\right)^2$$

we get the identity

(21) $$\frac{\left(1-tzC\left(t^2,z^2\right)\right)^2}{(1-tz)^2-z^2}=C\left(t^2,z^2\right),$$

which will be needed later.

From (20) we get

$$C\left(t,z^2\right)\left(1-z^2-tz^2C\left(t,z^2\right)\right)=1-tz^2C\left(t,z^2\right)$$

and thus

(22) $$C\left(t,z^2\right)=\frac{1-tz^2C\left(t,z^2\right)}{1-tz^2C\left(t,z^2\right)-z^2}=\frac{1}{1-\dfrac{z^2}{1-tz^2C\left(t,z^2\right)}}$$

This gives

**Example 4**

*For* $C\left(t,z^2\right)$ *we get* $s_n=0$ *and* $t_{2n}=1$ *and* $t_{2n+1}=t.$

*Thus* $C_n(t)$ *can be interpreted as the weight of Dyck paths of semi-length* $n$ *where the down-steps which end on even heights have weight* $1$ *and the other down-steps have weight* $t$*.*

**Example 5**

*For*

(23) $$G(t,z)=\sum_{n\geq 0}C_{n+1}(t)z^n=\frac{C(t,z)-1}{z}$$

*we get* $s_n=1+t$ *and* $t_n=t.$

*This also implies Coker's identity*

(24) $$C_{n+1}(t)=\sum_{k=0}^{\left\lfloor \frac{n}{2}\right\rfloor}\binom{n}{2k}C_k t^k(1+t)^{n-2k}.$$

**Proof**

Setting $C(t,z)=1+zG(t,z)$ in (19) gives

(25) $$G(t,z)=1+(1+t)zG(t,z)+tz^2G(t,z)^2$$

and therefore

(26) $$G(t,z)=\frac{1}{1-(1+t)z-tz^2G(t,z)}$$

which gives the continued fraction with $s_n=1+t$ and $t_n=t.$

Solving for $G(t,z)$ gives

(27) $$G(t,z)=\frac{1-(1+t)z-\sqrt{\left(1-(1+t)z\right)^2-4tz^2}}{2tz^2}.$$

In the same way as above we get Coker's identity (24).

For later use let us compute $\left(1-tz^2G\left(t^2,z^2\right)\right)^2$. By (25) we get

$$1-2tz^2G\left(t^2,z^2\right)+t^2z^4G\left(t^2,z^2\right)^2=-2tz^2G\left(t^2,z^2\right)+G\left(t^2,z^2\right)-\left(1+t^2\right)z^2G\left(t^2,z^2\right)$$
$$=G\left(t^2,z^2\right)\left(1-2tz^2-z^2-t^2z^2\right)$$

and therefore

$$\frac{\left(1-tz^2G\left(t^2,z^2\right)\right)^2}{1-(1+t)^2z^2}=G\left(t^2,z^2\right). \tag{28}$$

**Example 6**

*For* $C(t,z)$ *we get*

$s_n=1+t,\ \ s_0=1$ *and* $t_n=t.$

**Proof**

By (19) we get

$$C(t,z)=\frac{1}{1-(1-t)z-tzC(t,z)}=\frac{1}{1-(1-t)z-tz\left(1+zG(t.z)\right)}=\frac{1}{1-z-tz^2G(t.z)}$$

## 5. Analogous results for q-analogues at q=-1.

The most natural $q-$analogue of $C_n$ is $C_n(q)=\frac{1}{\left[n+1\right]}\begin{bmatrix}2n\\ n\end{bmatrix}$, where as usual

$\left[n\right]=\left[n\right]_q=\sum_{j=0}^{n-1}q^j,\ [n]!=\prod_{j=1}^{n}[j]$ and $\begin{bmatrix}n\\ k\end{bmatrix}=\frac{[n]!}{[k]![n-k]!}.$

In this case there are no nice formulas for $s_n$ and $t_n$. But there are nice ones if $q=-1,$ where we get $C_n(-1)=\binom{n}{\left\lfloor \frac{n}{2}\right\rfloor}$. We shall see in Theorem 3 that $\sum_{n\geq 0}C_n(-1)z^{2n}$ has a continued fraction with $s_n=0$ and $t_n=(-1)^{\binom{n}{2}}$.

We consider the $q-$Narayana polynomials

$$C_n(t;q)=\sum_{k=0}^{n-1}q^{k^2+k}\begin{bmatrix}n\\ k\end{bmatrix}\begin{bmatrix}n-1\\ k\end{bmatrix}\frac{t^k}{\left[k+1\right]} \tag{29}$$

which satisfy $C_n(1;q)=C_n(q)$ (cf. [6]).

Let

$$c_n(t) = C_n(t;-1) \tag{30}$$

be the limits of $C_n(t;q)$ for $q \to -1$.

The first terms are

$$\left(c_n(t)\right)_{n\geq 0} = \left(1,1,1+t,1+t+t^2,1+2t+2t^2+t^3,1+2t+4t^2+2t^3+t^4,\cdots\right).$$

As shown in [2] these polynomials are given by

$$c_n(t) = \sum_{k=0}^{n} v(n,k) t^k \tag{31}$$

with

$$v(n,k) = \binom{\left\lfloor \frac{n-1}{2} \right\rfloor}{\left\lfloor \frac{k}{2} \right\rfloor} \binom{\left\lfloor \frac{n}{2} \right\rfloor}{\left\lfloor \frac{k+1}{2} \right\rfloor} \tag{32}$$

and satisfy the recursion

$$\begin{aligned} c_{2n}(t) &= (1+t)c_{2n-1}(t), \\ c_{2n+1}(t) &= (1+t)c_{2n}(t) - tC_n\left(t^2\right). \end{aligned} \tag{33}$$

As analogue of $G(t,z)$ we consider

$$g(t,z) = \sum_{n\geq 0} c_{n+1}(t) z^n. \tag{34}$$

By (33) we get

$$\left(1-(1+t)z\right) g(t,z) = 1 - z^2 t G\left(t^2, z^2\right). \tag{35}$$

This implies

$$g(t,z) g(t,-z) = G\left(t^2, z^2\right) \tag{36}$$

because by (28)

$$g(t,z)g(t,-z) = \frac{\left(1-z^2tG\left(t^2,z^2\right)\right)}{1-(1+t)z} \frac{\left(1-z^2tG\left(t^2,z^2\right)\right)}{1+(1+t)z} = G\left(t^2,z^2\right).$$

Thus (35) and (36) give

$$\left(1-(1+t)z\right) g(t,z) = 1 - z^2 t g(t,z) g(t,-z) \tag{37}$$

which implies

(38) $$g(t,z)=\frac{1}{1-(1+t)z+tz^2 g(t,-z)}.$$

This gives

**Theorem 1**

*The generating function*

(39) $$g(t,z)=\sum_n c_{n+1}(t)z^n$$

*has a continued fraction expansion with* $s_n=(-1)^n\left(1+t\right)$ *and* $t_n=-t.$

*This also implies*

(40)
$$c_{2n+1}(t)=\sum_{k=0}^{n}(-t)^k\binom{n}{k}C_k(1+t)^{2n-2k},$$
$$c_{2n+2}(t)=\sum_{k=0}^{n}(-t)^k\binom{n}{k}C_k(1+t)^{2n-2k+1}.$$

The proof of (40) follows from the

**Lemma**

*If* $s_n=(-1)^n s$ *and* $t_n=t$ *then the weight of the Motzkin paths of length* $2n+\delta,\ \delta\in\{0,1\}$, *is given by* $\sum_{k=0}^{n}\binom{n}{k}C_k t^k s^{2n+\delta-2k}$.

The following proof is due to Martin Rubey [6].

**Proof.**

Consider a path of length $2n$ as a sequence of $n$ pairs of steps $U,H,D.$

Consider first the paths which do not have pairs with precisely one $H$. They consist of pairs $(U,U)$, $(U,D)$, $(D,U)$ and $(H,H)$. If the path has $k$ up-steps then there are $k$ pairs not of the form $(H,H)$. If we remove those there remains a Dyck path of semi-length $k.$ To obtain the weight of these paths there are $\binom{n}{k}$ choices for the pairs and $C_k$ Dyck paths on this set with weight $t^k$, which gives the formula. On the remaining paths we define a fixed-point free, sign-reversing involution. Find the first pair which contains precisely one $H$. Swapping this H with the step in the same pair changes the height of the H by one.

**Theorem 2**

*The continued fraction of*

*(41)* $$c(t,z)=\sum_{n\geq 0} c_n(t) z^n$$

*has* $s_0=1,\ s_n=(-1)^n\left(1-t\right)$ *and* $t_n=t.$

**Proof**

By (33) we get $\left(1-(1+t)z\right)c(t,z)=1-tzC\left(t^2,z^2\right).$

This implies by (21)

(42) $$c(t,z)c(-t,-z)=\frac{\left(1-tzC\left(t^2,z^2\right)\right)}{\left(1-(1+t)z\right)}\frac{\left(1-tzC\left(t^2,z^2\right)\right)}{\left(1+(1-t)z\right)}=C\left(t^2,z^2\right).$$

Therefore, we get

$\left(1-(1+t)z\right)c(t,z)=1-tzC\left(t^2,z^2\right)=1-tzc(t,z)c(-t,-z),$
which implies

(43) $$c(t,z)\left(1-(1+t)z+tzc(-t,-z)\right)=1.$$

Thus

$$c(t,z)=\frac{1}{1-(1+t)z+tzc(-t,-z)}=\frac{1}{1-z-tz^2g(-t,-z)}$$

which implies Theorem 2.

**Theorem 3**

The continued fraction of

(44) $$c(t,z^2)=\sum_{n\geq 0} c_n(t) z^{2n}$$

has $s_n=0,\ \ t_{2n}=(-1)^n, t_{2n+1}=(-1)^n t.$

**Proof**

By (43) we get

(45) $$c(t,z)\left(1-z-tzc(t,-z)\right)=1-tzc(t,-z)$$

which gives

(46) $$c(t,z^2)=\cfrac{1}{1-\cfrac{z^2}{1-tz^2c(t,-z^2)}}$$

As special case we get that

(47) $$c\left(z^2\right)=\sum_{n\geq 0} C_n(-1)z^{2n}$$

has a continued fraction expansion with $s_n=0$ and $t_n=\left(-1\right)^{\binom{n}{2}}$.

Let us finally mention that (46) implies

$$c(t,z)=\cfrac{1}{1-\cfrac{z}{1-tzc(t,-z)}}=\cfrac{1}{1-\cfrac{z}{1-\cfrac{tz}{1+\cfrac{z}{1+tzc(t,z)}}}}$$

which gives the explicit formula

(48) $$c(t,z)=\frac{\left(1+(1-t)z\right)}{2tz}\left(\sqrt{\frac{\left(1-(1-t)z\right)\left(1+(1+t)z\right)}{\left(1+(1-t)z\right)\left(1-(1+t)z\right)}}-1\right).$$

**6. Closed formulae**

Formula (48) sheds new light to some of the above and related results. If we set

(49) $$\gamma(t,z)=1+(1+t)z$$

it can be written as

(50) $$c(t,z)=\frac{\gamma(-t,z)}{2tz}\left(\sqrt{\frac{\gamma(t,z)\gamma(-t,-z)}{\gamma(t,-z)\gamma(-t,z)}}-1\right).$$

From this point of view formula (42) can also be obtained in the following way:

$$c(t,z)c(-t,-z)=\frac{\gamma(-t,z)}{2tz}\left(\sqrt{\frac{\gamma(t,z)\gamma(-t,-z)}{\gamma(t,-z)\gamma(-t,z)}}-1\right)\frac{\gamma(t,-z)}{2tz}\left(\sqrt{\frac{\gamma(t,z)\gamma(-t,-z)}{\gamma(t,-z)\gamma(-t,z)}}-1\right)$$

$$=\frac{\gamma(-t,z)\gamma(t,-z)}{4t^2z^2}\left(\sqrt{\frac{\gamma(t,z)\gamma(-t,-z)}{\gamma(t,-z)\gamma(-t,z)}}-1\right)^2$$

$$=\frac{1}{4t^2z^2}\left(\gamma(t,z)\gamma(-t,-z)+\gamma(-t,z)\gamma(t,-z)-2\sqrt{\gamma(t,z)\gamma(-t,-z)\gamma(t,-z)\gamma(-t,z)}\right).$$

Since $\gamma(t,z)\gamma(-t,-z)+\gamma(-t,z)\gamma(t,-z)=2\left(1-\left(1-t^2\right)z^2\right)$ and

(51) $$\begin{aligned}&\gamma(t,z)\gamma(-t,-z)\gamma(t,-z)\gamma(-t,z)\\&=\left(1-\left(1+t^2\right)z^2\right)^2-4t^2z^4=\left(1-\left(1-t^2\right)z^2\right)^2-4t^2z^2\end{aligned}$$

we get by (17)

$$c(t,z)c(-t,-z)=\frac{1}{2t^2z^2}\left(1-\left(1-t^2\right)z^2-\sqrt{\left(1-\left(1-t^2\right)z^2\right)^2-4t^2z^2}\right)=C\left(t^2,z^1\right).$$

It turns out that all square roots appearing in the above generating functions are square roots of rational functions where numerator and denominator are products of some of the functions $\gamma\left(t,z\right),\gamma\left(t,-z\right),\gamma\left(-t,z\right),\gamma\left(-t,-z\right)$ such that each of these functions occurs exactly once.

**Example A**

(52) $$C\left(t^2,z^2\right)=\frac{1}{2t^2z^2}\left(1-\left(1-t^2\right)z^2-\sqrt{\gamma(t,z)\gamma(-t,-z)\gamma(t,-z)\gamma(-t,z)}\right).$$

**Example a**

(53) $$c(t,z)=\frac{1}{2tz}\left(\sqrt{\frac{\gamma(t,z)\gamma(-t,-z)\gamma(-t,z)}{\gamma(t,-z)}}-\gamma(-t,z)\right).$$

Let $W_n(t)=\sum_{k=0}^{n}\binom{n}{k}^2t^k$ and (cf. e.g. [5])

(54) $$W(t,z)=\sum_{n\geq 0}W_n(t)z^n==\frac{1}{\sqrt{1-2(t+1)z+(t-1)^2z^2}}.$$

By (51) we get

**Example B**

(55) $$W\left(t^2,z^2\right)=\frac{1}{\sqrt{\gamma(t,z)\gamma(-t,-z)\gamma(t,-z)\gamma(-t,z)}}.$$

Let $W_n(t;q)=\sum_{k=0}^{n}\begin{bmatrix}n\\k\end{bmatrix}^2t^k$ and $w_n(t)=W_n\left(t;-1\right)$ the value at $q=-1.$ It can easily be verified that

(56) $$\begin{aligned}w_{2n}(t)&=W_n\left(t^2\right),\\ w_{2n+1}(t)&=(1+t)w_{2n}(t).\end{aligned}$$

This implies that $w(t,z)=\sum_{n\geq 0}w_n(t)z^n$ satisfies

(57) $$w(t,z)=\gamma(t,z)W\left(t^2,z^2\right)=\frac{\gamma(t,z)}{\sqrt{\gamma(t,z)\gamma(-t,-z)\gamma(t,-z)\gamma(-t,z)}}.$$

Thus we get

**Example b**

(58)
$$w(t,z)=\sqrt{\frac{\gamma(t,z)}{\gamma(-t,z)\gamma(t,-z)\gamma(-t,-z)}}.$$

By (33) the numbers $c_n(t)$ are uniquely determined by the numbers $c_{2n+1}(t)$ which by [1](13) satisfy $c_{2n+1}(t)=W_n\left(t^2\right)+ntC_n\left(t^2\right).$

Let $h(t,z)=\sum_{n\geq 0}c_{2n+1}(t)z^n.$ Here we get

**Example c**

(59)
$$h(t,z^2)=\frac{1}{2tz^2}\left(\sqrt{\frac{\gamma(-t,z)\gamma(-t,-z)}{\gamma(t,z)\gamma(t,-z)}}-1\right).$$

**Proof**

By (33) we have

$\gamma(t,z)zh(t,z^2)=\gamma(t,z)\sum_{n\geq 0}c_{2n+1}(t)z^{2n+1}=\sum_{n\geq 0}c_{n+1}(t)z^{n+1}=c(t,z)-1$ which implies

$$h(t,z^2)=\frac{c(t,z)-1}{z\gamma(t,z)}=\frac{1}{2tz^2\gamma(t,z)}\left(\sqrt{\frac{\gamma(t,z)\gamma(-t,-z)\gamma(-t,z)}{\gamma(t,-z)}}-\gamma(-t,z)-2tz\right)$$
$$=\frac{1}{2tz^2}\left(\sqrt{\frac{\gamma(-t,-z)\gamma(-t,z)}{\gamma(t,z)\gamma(t,-z)}}-1\right)=\frac{1}{2tz^2}\left(\sqrt{\frac{1-(1-t)^2z^2}{1-(1+t)^2z^2}}-1\right).$$

This implies by Example A

$$\frac{1}{h\left(t,z^2\right)}=2tz^2\left(\sqrt{\frac{\gamma(-t,-z)\gamma(-t,z)}{\gamma(t,z)\gamma(t,-z)}}+1\right)\frac{\gamma(t,z)\gamma(t,-z)}{4tz^2}=$$
$$\frac{1}{2}\left(\sqrt{\gamma(-t,-z)\gamma(-t,z)\gamma(t,z)\gamma(t,-z)}+1-(1+t)^2z^2\right)=1-(1+t)z^2-t^2z^2C\left(t^2,z^2\right)$$
$$1-\left(1+t+t^2\right)z^2-t^2z^4G\left(t^2,z^2\right)$$

and therefore, $\frac{1}{h(t,z)}=1-\left(1+t+t^2\right)z-t^2z^2G\left(t^2,z\right).$

By Example 5 we get

**Theorem 4**

The generating function

(60)
$$h(t,z)=\sum_{n\ge 0} c_{2n+1}(t)z^n$$

has a continued fraction expansion with

(61)
$$s_0=1+t+t^2,\ \ s_n=1+t^2 \text{ for } n>0,$$
$$t_n=t^2.$$

## References

[1] Johann Cigler, Symmetric Dyck paths and q-Narayana numbers, arXiv:2601.08366

[2] Johann Cigler, Some remarks about q-Narayana polynomials for q=-1, arXiv:2603.04137

[3] P. Flajolet, Combinatorial aspects of continued fractions, Discr. Math. 32(1980),125-161

[4] Christian Krattenthaler, Advanced Determinant Calculus, Séminaire Lotharingien Combin. 42 (1999), B42q

[5] Qiongqiong Pan and Jiang Zeng, On total positivity of Catalan-Stieltjes matrices. Electronic J. Comb. 23(4), 2016, P4.33

[6] T. Kyle Petersen, Eulerian numbers, Birkhäuser 2015

[7] Martin Rubey, Answer to https://mathoverflow.net/questions/509835/